# Nishida Relations in Bordism and Homology
(June 1995. To Appear, C.R.Math.Rep.Acad.Sci.Canada))

Terrence Bisson     André Joyal
bisson@canisius.edu     joyal@math.uqam.ca

This is the second of a series of Compte Rendus. In the first [1] we have presented a theory of Dyer-Lashof operations in unoriented bordism. Here we shall discuss the (Nishida) relations between Dyer-Lashof and Landweber-Novikov operations. They are used to represent the algebra $N_*\Sigma$ of covering manifolds in terms of their homology characteristic numbers. The proofs are based on the properties of the covering space operations and the notions of $D$-ring and $Q$-ring introduced in [1].

## 1. The Nishida relations in homology.

In homology mod 2 the Nishida relations are commutation relations between Dyer-Lashof and Steenrod operations (see [4] for instance). We shall express the Nishida relations as a commutative square which involves the Milnor coaction and the $Q$-structure on the homology of any $E_\infty$-space discussed in [1].

Recall that the Milnor Hopf algebra $\mathcal{A}_*$ is the dual of the Steenrod algebra; see [6] for instance. As a graded algebra it is $\mathcal{A}_* = \mathbf{Z}_2[\xi_0^\pm, \xi_1, \xi_2, \ldots] = \xi_0^{-1}\mathbf{Z}_2[\xi_0, \xi_1, \xi_2, \ldots]$ with $\operatorname{grade}(\xi_i) = 2^i - 1$; the diagonal $\delta : \mathcal{A}_* \to \mathcal{A}_* \otimes \mathcal{A}_*$ is the unique ring homomorphism such that $\delta(\xi) = (\xi \otimes 1) \circ (1 \otimes \xi)$ where $\xi(x) = \sum \xi_i x^{2^i}$. We are diverging from the usual convention that puts $\xi_0 = 1$.

The homology of any space has a natural (left) coaction

$$\alpha : H_*(X) \to \mathcal{A}_* \otimes H_*(X) = H_*(X)[\xi_0^\pm, \xi_1, \xi_2, \ldots]$$

which restricts to the usual (left) coaction when we put $\xi_0 = 1$ and to the "grading" coaction $\alpha(x_n) = \xi_0^{-n} x_n$ for $x \in H_n(X)$ when we put $0 = \xi_1 = \xi_2 = \cdots$.

For example, for $X = RP^\infty$ we have $\alpha(b)(x) = b(\xi^{(-1)}(x))$ where $b(x) = \sum_i b_i x^i$ and $b_0, b_1, \ldots$ is the canonical basis of $H_*RP^\infty$ and $\xi^{(-1)}(x)$ is the composition inverse of the power series $\xi(x)$.

**Proposition.** *If $R$ is a $Q$-ring, then there is a $Q$-structure on $\mathcal{A}_* \otimes R$ determined by $Q_t(\xi)(x(x+t)) = \xi(x)\xi(x+t)$. In particular, there is a unique $Q$-structure on $\mathcal{A}_*$ such that $Q_t(\xi)(x(x+t)) = \xi(x)\xi(x+t)$.*

**Remark:** Since the definition forces

$$Q_t(\xi_0) = \xi_0 \sum_i \xi_i t^{2^i - 1}$$

we see that the $Q$-structure on $\mathcal{A}_*$ cannot survive if we try to insist that $\xi_0 = 1$.

Suppose now that $R$ is both a $Q$-ring and has a Milnor coaction. We give $\mathcal{A}_* \otimes R$ the $Q$-structure from the proposition.

**Definition:** The *Nishida relations* hold for a $Q$-ring $R$ if the following diagram

$$\begin{array}{ccc} R & \xrightarrow{Q_t} & R[[t]] \\ \alpha \downarrow & & \downarrow \alpha \\ \mathcal{A}_* \otimes R & \xrightarrow{Q_t} & \mathcal{A}_* \otimes R[[t]] \end{array}$$

commutes, where $\alpha$ is extended to $R[[t]]$ by putting $\alpha(t) = \xi(t)$.

**Theorem.** *If $X$ is an $E_\infty$-space, then $H_*(X)$ is a $Q$-ring with a Milnor coaction, and the Nishida relations hold for $H_*(X)$.*

This is a complete description of the Nishida relations in mod 2 homology. We want to give a similar description in unoriented bordism

## 2. The Nishida relations in unoriented bordism.

For any $\mathbf{Z}_2$-algebra $R$ let $\mathcal{B}_*(R)$ be the group of inversible formal power series $f(x) \in xR[[x]]$ under substitution. The functor $R \mapsto \mathcal{B}_*(R)$ is representable by a Hopf algebra $\mathcal{B}_*$ called the *Faa di Bruno Hopf algebra* (see example 49 in [2] for instance). We have $\mathcal{B}_* = \mathbf{Z}_2[h_0^\pm, h_1, \ldots] = h_0^{-1}\mathbf{Z}_2[h_0, h_1, \ldots]$ and the diagonal of $\mathcal{B}_*$ is given by

$$\delta(h) = (h \otimes 1) \circ (1 \otimes h)$$

where $h(x) = \sum_n h_n x^{n+1}$. We want to consider left comodules over $\mathcal{B}_*$. Recall that the tensor product of left comodules over a bialgebra is a left comodule.

**Example:** Let $N_*$ denote the Lazard ring for formal power series of order two (see [5] for example). There is a natural (left) coaction $\phi : N_* \to \mathcal{B}_* \otimes N_*$ encoding change of parameters in formal group laws. We have $\phi(F)(x,y) = h(F(h^{-1}x, h^{-1}y))$ where $\phi : N_* \to N_*[h_0^\pm, h_1, \ldots]$.

Notice that $N_* \otimes N_* \to N_*$ is a morphism of left $\mathcal{B}_*$ comodules, so that $N_*$ is a monoid in the category of left $\mathcal{B}_*$ comodules.

**Definition:** A *Landweber-Novikov coaction* is a module over $N_*$ in the category of $\mathcal{B}_*$-comodules. More explicitly, this is a module over $N_*$ with a comodule structure over $\mathcal{B}_*$ such that the module structure map $N_* \otimes M \to M$ is a map of $\mathcal{B}_*$-comodules.

In order to describe the natural Landweber-Novikov coaction on unoriented bordism, we first need to recall the theory of cobordism characteristic classes for vector bundles, as sketched by Quillen in [5] for instance; we use a slightly unstable version of the usual definitions, however. For any space $X$ the total characteristic class of a real vector bundle $V$ on $X$ is the element $c(V) \in N^*(X)[b_0, b_1, ...]$ having the following properties:
  1) the map $V \to c(V)$ is a natural transformation $c : Vect(X) \to N^*(X)[b_0, b_1, ...]$
  2) $c(V_1 \oplus V_2) = c(V_1)c(V_2)$

3) $c(L) = \sum b_i e(L)^i$ if $L$ is a line bundle with euler class $e(L) \in N^1(X)$. We have an expansion
$$c(V) = \sum_R c_R(V) b^R$$
where $b^R = b_0^{r_0} b_1^{r_1} \cdots$ for $R = (r_0, r_1, \ldots)$. For virtual bundles the total characteristic class is the element of $N^*(X)[b_0^\pm, b_1, \ldots]$ defined by putting $c(V - W) = c(V)c(W)^{-1}$.

**Example:** The unoriented bordism group $N_*(X)$ of any space $X$ has a Landweber-Novikov coaction given by the Landweber-Novikov (total) operation
$$\phi : N_*(X) \to \mathcal{B}_* \otimes N_*(X) = N_*(X)[h_0^\pm, h_1, \ldots].$$
It can be defined via the following explicit formula. If $f : M \to X$ represents an element of $N_*(X)$ then
$$\phi(M, f) = \sum_R f_*(c_R(\nu_M) \cap \mu_M) h^R$$
where $\nu_M = -\tau_M$ is the normal bundle of $M$ and $\mu_M$ is the fundamental class of $M$.

**Example:** The Landweber-Novikov coaction on $N_*(RP^\infty)$ is characterised by the identity $\phi(b)(x) = b(h^{(-1)}(x))$ where $b(x) = \sum_i b_i x^i$.

**Proposition.** If $R$ is a $D$-ring, then there is a $D$-structure on $\mathcal{B}_* \otimes R = R[h_0^\pm, h_1, \ldots]$ determined by $D_t(h)(x(x+t)) = h(x)h(F(x,t))$.

Suppose now that $R$ is a $D$-ring which also has a Landweber-Novikov coaction. We give $\mathcal{B}_* \otimes R$ the $D$-structure from the proposition.

**Definition:** The *Nishida relations* hold for a $D$-ring $R$ if the diagram

$$\begin{array}{ccc} R & \xrightarrow{D_t} & R[[t]] \\ \phi \downarrow & & \downarrow \phi \\ \mathcal{B}_* \otimes R & \xrightarrow{D_t} & \mathcal{B}_* \otimes R[[t]] \end{array}$$

commutes, where $\phi$ is extended to $R[[t]]$ by putting $\phi(t) = h(t)$.

Suppose that $X$ is an $E_\infty$-space. Then $N_* X$ is an $D$-ring together with a Landweber-Novikov coaction.

**Theorem 2.** *If $X$ is an $E_\infty$-space, then the Nishida relations hold for $N_* X$ in that the following diagram commutes:*

$$\begin{array}{ccc} N_*(X) & \xrightarrow{D_t} & N_*(X)[[t]] \\ \phi \downarrow & & \downarrow \phi \\ \mathcal{B}_* \otimes N_*(X) & \xrightarrow{D_t} & \mathcal{B}_* \otimes N_*(X)[[t]] \end{array}$$

Here $\phi$ has been extended to $N_*(X)[[t]]$ by putting $\phi(t) = h(t)$.

Recall from [1] that for any space $X$ we have $N_*E_\infty(X) = D\langle N_*(X)\rangle$, the free $D$-ring generated by $N_*(X)$. Of course, $N_*(X)$ also has a Landweber-Novikov coaction.

**Theorem 3.** *Suppose that $M$ has a Landweber-Novikov coaction, and let $D\langle M\rangle$ be the free $D$-ring generated by $M$. Then there is a unique Landweber-Novikov coaction on $D\langle M\rangle$ which satisfies all of the following:*
  i) *the canonical map $M \to D\langle M\rangle$ is a comodule map;*
  ii) *$D\langle M\rangle$ is an algebra over its Landweber-Novikov coaction;*
  iii) *$D\langle M\rangle$ satisfies the Nishida relations.*

The theorem shows that the Landweber-Novikov coaction on $N_*E_\infty(X)$ is fully determined by its values on $N_*X$. In particular, the Landweber-Novikov coaction on $N_*(B\Sigma_*) = D\langle x\rangle$ is determined by the relation $\phi(x) = x$.

It is well known that the Thom reduction $\epsilon : N_*(X) \to H_*(X)$ induces an isomorphism
$$N_*(X) \otimes_{N_*} \mathbf{Z}_2 \simeq H_*(X)$$

We shall write $\epsilon : \mathcal{B}_* \to \mathcal{A}_*$ for the ring homomorphism such that $\epsilon(h) = \xi$. We will characterize the behavior of the Thom reduction with respect to the usual Landweber-Novikov operations in unoriented cobordism and the Steenrod operations in mod 2 homology by showing that the Thom reduction gives rise to a simple equivalence of categories.

Let $[\mathcal{B}_\star \mathcal{N}_\star]$ denote the category of Landweber-Novikov coactions (left $\mathcal{B}_*$-comodules which are modules over $N_*$) and let $[\mathcal{A}_\star]$ denote the category of Milnor coactions (left $\mathcal{A}_*$-comodules).

For any $M \in [\mathcal{B}_\star \mathcal{N}_\star]$ let us put $T(M) = M \otimes_{N_*} \mathbf{Z}_2$. By composing the coaction $M \to \mathcal{B}_* \otimes M$ with $\epsilon : \mathcal{B}_* \to \mathcal{A}_*$ we obtain a coaction $M \to \mathcal{A}_* \otimes M$. By further tensoring with $\mathbf{Z}_2$ we obtain a coaction $T(M) \to \mathcal{A}_* \otimes T(M)$. This defines a functor $T : [\mathcal{B}_\star \mathcal{N}_\star] \to [\mathcal{A}_\star]$.

**Proposition.** *The functor $T$ defines an equivalence of categories $[\mathcal{B}_\star \mathcal{N}_\star] \simeq [\mathcal{A}_\star]$.*

In fact, we get a similar result even when we take the Dyer-Lashof operations into account, which we do by using monads which encode the Dyer-Lashof operations in bordism and homology. For background on monads and their category of actions see [3] for instance.

The functor $M \mapsto D\langle M\rangle$ from $[\mathcal{B}_\star \mathcal{N}_\star]$ to itself is a monad that we shall denote $D$ (it is a monad since $D\langle M\rangle$ is a free structure). The $D$-actions are exactly the $D$-rings with a Landweber-Novikov coaction which satisfies the Nishida relations. Let us denote this category by $[\mathcal{B}_\star \mathcal{N}_\star]^D$. Similarly, the functor $M \mapsto Q\langle M\rangle$ from $[\mathcal{A}_\star]$ to itself is a monad that we shall denote $Q$. The $Q$-actions are exactly the $Q$-rings with a Milnor coaction which satisfies the Nishida relations. Let us denote this category by $[\mathcal{A}_\star]^Q$.

**Proposition 2.** *The functor $T$ transforms the $D$-actions into the $Q$-actions and defines an equivalence of categories*
$$[\mathcal{B}_\star \mathcal{N}_\star]^D \simeq [\mathcal{A}_\star]^Q.$$

This result shows that for any $E_\infty$-space $X$, the $D$-ring $N_*(X)$ can be recovered entirely from the $Q$-ring $H_*(X)$ as long as the coaction of $\mathcal{A}_*$ on $H_*(X)$ is known.

## 3. Characteristic numbers and covering space operations

The bordism classes of manifolds are determined by their tangential characteristic numbers, or equally by their normal characteristic numbers.

Recall the discussion in section 2 of characteristic classes in unoriented cobordism. The total Stiefel-Whitney class $w(V) \in H^*(X)[b_0, b_1, ...]$ of a vector bundle $V$ on $X$ is the Thom reduction of the total cobordism characteristic class. Then $w(\ )$ is multiplicative in that $w(V_1 \oplus V_2) = w(V_1)w(V_2)$, and $w(L) = \sum_i b_i e(L)^i$ for any line bundle over $X$, where $e(L) \in H^1(X)$ is the euler class of $L$.

The tangential characteristic numbers can be grouped together as the coefficients of a polynomial $\beta^\tau(M) \in \mathbf{Z}_2[h_0, h_1, ...]$. We have

$$\beta^\tau(M) = \sum_R \langle w_R(\tau_M), \mu_M \rangle h^R$$

where $w(\tau_M) = \sum_R w_R(\tau_M) b^R$ is the total Stiefel-Whitney class of the tangent bundle of $M$ and $\mu_M$ is the fundamental class of $M$.

In keeping with the Landweber-Novikov coaction, we can work instead with characteristic number polynomials for the normal bundle, the virtual bundle $\nu_M = -\tau_M$. More generally, as in [5] for instance, we can define a *Boardman map*

$$\beta: N_*(X) \to \mathcal{B}_* \otimes H_*(X) = H_*(X)[h_0^\pm, h_1, ...]$$

by the explicit formula

$$\beta(M, f) = \sum_R f_*(w_R(\nu_M) \cap \mu_M) h^R.$$

Let $B\Sigma_*$ denote the classifying space for finite covering spaces. It is an $E_\infty$-space. Let $N_*\Sigma$ denote $N_*B\Sigma_*$ and $H_*\Sigma$ denote $H_*B\Sigma$. The Boardman map $\beta: N_*\Sigma \to \mathcal{B}_* \otimes H_*\Sigma$ obviously preserves sums and products. In fact, it also preserves the operation of substitution, when substitution is properly defined on $\mathcal{B}_* \otimes H_*\Sigma$.

**Proposition.** *If $R$ is a $Q$-ring then there is a $Q$-structure on $\mathcal{B}_* \otimes R$ determined by $Q_t(h)(x(x+t)) = h(x)h(x+t)$. In particular, there is a unique $Q$-structure on $\mathcal{B}_*$ such that $Q_t(h)(x(x+t)) = h(x)h(x+t)$.*

From [1] we have that $H_*\Sigma = Q\langle x \rangle$ is the free $Q$-ring on one generator. Since $\mathcal{B}_* \otimes H_*\Sigma$ is the coproduct in the category of $Q$-rings, we can view $\mathcal{B}_* \otimes H_*\Sigma$ as the collection of unary operations in the theory of $Q$-rings which are extensions of $\mathcal{B}_*$. Therefore $H_*\Sigma[h_0^\pm, h_1, ...] = \mathcal{B}_* \otimes H_*\Sigma$ admits an operation of substitution. Here is an explicit formula for substituting two elements of $H_*\Sigma[h_0^\pm, h_1, ...]$:

$$\left(\sum_R p_R h^R\right) \circ \left(\sum_R q_R h^R\right) = \sum_{R,S} p_R(q_S h^S) h^R$$

where the operations $p_R$ are applied to the elements $q_S h^S$ of the $Q$-ring $\mathcal{B}_* \otimes H_*\Sigma$.

**Theorem 4.** *The Boardman map $\beta : N_*\Sigma \to H_*\Sigma[h_0^{\pm}, h_1, ...]$ preserves sums, products and substitutions.*

We view this as describing for each covering $p$ of closed manifolds how the normal characteristic numbers of $p(M)$ are determined by the characteristic numbers of $M$. The complete description is summarized by the formula

$$\beta(p(M)) = (\beta(p))(\beta(M)),$$

where the right-hand substitution takes place in $H_*\Sigma[h_0^{\pm}, h_1, \ldots]$.

**References**

[1] T. P. Bisson, A. Joyal, The Dyer-Lashof algebra in bordism, C. R. Math. Rep. Acad. Sci. Canada (to appear).
[2] A. Joyal, Une théorie combinatoire des séries formelles, Adv. in Math. 42 (1981), 1-82.
[3] S. MacLane. Categories for the Working Mathematician, Springer 1971.
[4] J. P. May, Homology operations on infinite loop spaces, in Proc. Symp. Pure Math. 22, A. M. S. 1971, 171-185.
[5] D. Quillen, Elementary proofs of some results of cobordism theory using Steenrod operations, Adv. in Math. 7 (1971), 29-56.
[6] N. E. Steenrod, D. B. A. Epstein, Cohomology operations. Ann. of Math. Studies no. 50, Princeton 1962.

(*) Canisius College, Buffalo, N.Y. (U.S.A). e-mail: *bisson@canisius.edu*.

(**) Département de Mathématiques, Université du Québec à Montréal, Montréal, Québec H3C 3P8. e-mail: *joyal@math.uqam.ca*.